\documentclass[11pt]{article}
\usepackage{graphicx}
\usepackage{amsmath}
\usepackage{amssymb}
\usepackage{amsthm}
\usepackage{epsfig}
\usepackage{latexsym} 
\usepackage{color}

\newtheorem{theorem}{Theorem}[section]
\newtheorem{lemma}[theorem]{Lemma}

\newtheorem{definition}[theorem]{Definition}

\title{Incompressible surfaces and spunnormal form} 
\author{G. S. Walsh \footnote{Supported in part by N. S. F. grant 0805908  }} 
\date{} 
\begin{document} 

\maketitle
\abstract{Suppose $M$ is a cusped finite-volume hyperbolic 3-manifold and  $\mathcal{T}$ is an ideal triangulation of $M$ with essential edges. We show that any incompressible surface $S$ in $M$ that is not a virtual fiber can be isotoped  into spunnormal form in $\mathcal{T}$.  The proof is based directly on ideas of W. Thurston. }
\vskip .2 in

\section{Definitions and background}

It is well known that an incompressible surface in a closed 3-manifold can be isotoped to be in normal form with respect to a triangulation of that 3-manifold.  Here we extend this result to show that an incompressible surface which is not a fiber or a semi-fiber in a hyperbolic 3-manifold with cusps can be isotoped to be in spunnormal form with respect to an ideal triangulation of that hyperbolic 3-manifold.  We require that the ideal triangulation contain only essential edges.  Ideal triangulations are extremely useful in the study of hyperbolic 3-manifolds.  For example, the program Snappea (now SnapPy) \cite{snappea, SnapPy} computes an ideal triangulation of a cusped hyperbolic 3-manifold and uses it to determine a wealth of data about the 3-manifold.  See also \cite{Choi, pflat}.

Here the proof does not use the hyperbolic structure, but only topological properties of the 3-manifold, namely that it is atoroidal, acylindrical and irreducible, and has torus boundary components. The interior of such a manifold admits a hyperbolic structure but we do not use this fact.  Thus we give the definitions for general 3-manifolds with boundary, and hope that this will not distract the reader.  

\begin{definition}  \label{ideal} Let $M$ be a 3-manifold with boundary.  An \emph{ideal triangulation} of $M$ is a map $\mathcal{T}$ from $M$ onto a pseudo-manifold $N$. A \emph{pseudo-manifold} is the quotient of a finite set of 3-simplices with 2-dimensional faces of 3-simplices identified in pairs by piecewise linear homeomorphisms. We require that that $\mathcal{T}$ restricts to a homeomorphism from the interior of $M$ to the complement of the 0-skeleton of $N$, and that each boundary component of $M$ maps to a vertex of $N$.  We further require that each edge have valence at least 3.

\end{definition} 
Note: This is slightly different than the usual definition of pseudo-manifold; we allow 2-dimensional faces of the same simplex to be identified.  The requirement on the valence of the edges is to ensure that the induced triangulation on the link of the vertex is non-degenerate.  An ideal triangulation naturally gives a PL-structure on the interior of $M$, and we work in the PL category. 

We denote a 3-manifold with an ideal triangulation by $(M, \mathcal{T})$, or sometimes $\mathcal{T}:M \rightarrow N$ for clarity.  A 3-simplex of a pseudo-manifold is the image of a single 3-simplex in the quotient map, and an ideal tetrahedron of $\mathcal{T}$ is the pre-image of a 3-simplex of $N$ in the interior of $M$.  Similarly, the ideal faces and edges of $\mathcal{T}$ are the pre-images of faces and edges of $N$.  We call the vertices of $N$ \emph{ideal vertices} of $(M, \mathcal{T})$. The boundary components of $M$ are in 1-1 correspondance with the ideal vertices. We will also place an \emph{end structure} on $(M, \mathcal{T})$, which is defined as follows.  There is a canonical neighborhood of each ideal vertex $v$ formed by taking the star of $v$ in the second barycentric subdivision of $N$.  The preimage of such a neighborhood is called an \emph{end}, and it has a natural product structure $L \times [0, \infty)$, where $L$ is the link of $v$ in the second barycentric subdivision of $N$.  In this natural product structure, vertical rays $x \times [0, \infty)$ are mapped linearly to rays which are linear in the star of $v$ and have $v$ as the deleted endpoint.  Label the ends of $(M, \mathcal{T})$ by $\mathcal{E}_i$, and call their union $\mathcal{E}$. For each $\mathcal{E}_i \neq S^2 \times \mathbb{R}_+$, fix an identification of the universal cover of $\mathcal{E}_i$ with the closed upper half space $\lbrace (x,y,z): z \geq 0 \rbrace$.  We require  that the universal cover of $L \times \lbrace 0 \rbrace$ is identified with the plane $z=0$, and that the horizontal planes $z =r$ are invariant under the action of $\pi_1(L)$.   This fixes a product structure on each end. Denote $ \overline{N \setminus \mathcal{E}}$ by $N_0$.  If $(M, \mathcal{T})$ is an ideal triangulation of $M$, then there is a homeomorphism $\phi_\mathcal{T}:N_0 \rightarrow M$.  


\begin{definition} An ideal triangulation $\mathcal{T} :M \rightarrow N$ has \emph{essential edges} if, for each edge $e$ of $N$, the proper arc $e \cap N_0$ is not homotopic rel endpoints into $\partial N_0$.  \end{definition} 

If the interior of $M$ admits a finite-volume hyperbolic structure then $M$ admits an ideal triangulation with ideal edges which are geodesics in this hyperbolic structure, although some of the tetrahedra may be degenerate, see \cite{EpPenn} and \cite{pflat}. This can be fattened up to make an ideal triangulation with essential edges. 

\begin{definition} An immersion $f: S \rightarrow M$ is a \emph{virtual fiber} if there is finite cover $p: \tilde M \rightarrow M$  and an embedded surface $\tilde S$ in $\tilde M$  such that $p(\tilde S) = f(S)$ and $\tilde S$ is a fiber in a fibration of $\tilde M$ over the circle. \end{definition} 

An embedded virtual fiber of $M$ is either a fiber of a fibration of $M$ over $S^1$, or a generic fiber of an orbifold-fibration of $M$ over an interval with mirrored endpoints. In this latter case, the surface is called a semi-fiber, and $M$ is the union of two twisted $I$-bundles over the surface. 

For the remainder of the paper, all surfaces in manifolds are embedded unless otherwise indicated. Thus virtual fibers are either fibers or semi-fibers. 

 In the following definition we consider faces which have been identified in the pseudo-manifold to be distinct.  
 
 \begin{definition}A properly embedded arc on the face of a 3-simplex is called \emph{normal} if it has endpoints on two different edges.  A \emph{normal curve} is a simple closed curve that is a union of normal arcs such that adjacent arcs lie on different faces. Note that if a disk in a tetrahedron has boundary a normal curve, then that curve has length 3 or 4, or length $\geq 8$. Disks with boundary a normal curve of length 3 or 4, respectively, in a tetrahedron are called \emph{normal triangles} or \emph{normal quadrilateral disks}, respectively. 
 \end{definition}
 
  We denote the interior of a surface $\Sigma$ by $\dot{\Sigma}$. 
 
 \begin{definition}Let $\Sigma$ be a surface, possibly with boundary.  The image of a map $g$ from $\Sigma$ into a pseudo-manifold $N$  is  \emph{spunnormal} if: 
 \begin{enumerate}
  \item The map $g$ restricted to $\dot{\Sigma}$ is an embedding.
  \item The preimage $g^{-1}(N^{(0)}) = \partial \Sigma$.
 \item The intersection of $g(\dot{\Sigma})$ with every tetrahedron of $N$ is a union of normal triangles and normal quadrilateral disks. We  require that $g(\dot{\Sigma})$ has finitely many quadrilateral disks.  
 \end{enumerate} 
 
 \end{definition}

We will often denote both the surface and its image by $\Sigma$.  Spunnormal surfaces were developed by W. Thurston as a natural extension of normal surfaces for manifolds with ideal triangulations.  This form is particularly applicable for surfaces in hyperbolic manifolds with cusps.   The space of spunnormal surfaces in a hyperbolic  3-manifold with an ideal triangulation can be easily computed from the edge equations for the triangulation. This follows from the fact that the surface must glue up around an edge. This is used for computations in a program by Culler and Dunfield  \cite{cullerprogram}.  The theory of the space of spunnormal surfaces has been developed in work  of Tillmann \cite{tillmann}, and Kang and Rubinstein in \cite{Rubspun}. 

A compact manifold $M$ is \emph{atoroidal} if every rank two abelian subgroup of $\pi_1(M)$ is conjugate into the image of the inclusion of the fundamental group of some boundary component, and $M$ is not an $I$-bundle over the torus or the Klein bottle, or a disk bundle over $S^1$.  $M$ is \emph{acylindrical} if it does not contain essential annuli.  A properly embedded two-sided surface $S$ in $M$ is \emph{incompressible} if, for each component $S_i$, either (1) $S_i$ is a disk with $\partial S_i$ an essential curve of $\partial M$ or (2) $\chi(S_i) \leq 0$ and $S_i$ does not admit any  compressing disks or boundary-compressing disks.  We say that an isotopy $F: S \times I \rightarrow M\times I$ is \emph{proper} if $F(S,t)$ is a proper map for each $t \in I$.  We are now able to state our main theorem, which is proven in Section \ref{finite}: 

\begin{theorem} \label{main}
Let $M$ be an atoroidal, acylindrical, irreducible, compact  3-manifold with torus boundary components, and $\mathcal{T}:M \rightarrow N$ an ideal triangulation of $M$ with essential edges.  Let $S$ be a properly embedded, two-sided, incompressible surface in $M$ that is not a virtual fiber.  Then there is a spunnormal surface $\Sigma$ in $N$ such that $\phi_{\mathcal{T}}(\Sigma \cap N_0) $ is properly isotopic to $S$. \end{theorem}

  Kneser proved that surfaces as above in a closed 3-manifold $M$ can be isotoped into normal form in any triangulation of $M$, and this proof holds for closed surfaces in $(M, \mathcal{T})$ as above. W. Thurston \cite{Thurstonpf} proved that 3-manifolds which satisfy the hypotheses are exactly the manifolds with torus boundary whose interiors admit finite-volume hyperbolic structures.  By work of F. Bonahon and W. Thurston, the surface subgroups without accidental parabolics that are not virtual fiber subgroups are exactly the quasi-Fuchsian subgroups. Therefore, Theorem \ref{main} implies that a quasi-Fuchsian surface in a finite volume hyperbolic 3-manifold with cusps can be realized in normal or spunnormal form.
   In \cite{Kangfig8}, it was shown that the condition that $S$ is not a fiber is necessary.   In particular, the fiber of the figure-8 knot complement cannot be realized in spunnormal form in the standard ideal triangulation.   Note that the hypotheses of the theorem imply that $M$ is $\partial$-irreducible (since $M$ is not a solid torus).   

\vskip .2 in 
\noindent
{\bf Acknowledgements} The author wishes to thank the anonymous referee for many criticisms and suggestions which improved the paper, and Nathan Dunfield for correcting Definition 1.2.

\section{Product regions} 
 For the purposes of this paper, the most important aspect of the assumption that $S$ is not a virtual fiber is that  product regions have bounded length with respect to $S$.  This idea has been used in several other contexts, namely  \cite{cooperlong1}, and \cite{li}, and in initial versions of the Cyclic Surgery Theorem, see \cite{shalenfirst}. 
 
 \label{product}
Let $M$ be a irreducible 3-manifold with boundary.  Let $\mathcal{A} $ be a collection of closed annuli in $\partial M$ such that  (1) $\overline{\partial M \setminus \mathcal{A}} $ is incompressible in $M$ and (2) essential arcs of $\mathcal{A} $ cannot be homotoped in $M$ rel boundary  into $\overline{\partial M \setminus \mathcal{A}} $.  

\begin{definition} \label{productdef} A \emph{product disk with respect to $\mathcal{A} $} is a proper map $i : I \times I \rightarrow M$ such that  $i(\partial I \times I)$ is a pair of essential arcs in $\mathcal{A} $ and $i(I \times \partial I)$ is a pair of possibly immersed arcs in $\overline{\partial M \setminus \mathcal{A}} $ that cannot be homotoped in $M$ rel boundary into $ \mathcal{A} $. 
\end{definition}

We call $i(\partial I \times I)$ the {\it vertical  boundary} of $D$ and  $i(I \times \partial I)$ the {\it horizontal boundary} of $D$.  All proper arcs isotopic to essential arcs of $\mathcal{A} $ are called {\it vertical arcs}.  We say that the image of $\rho: S \times I \rightarrow M$ is an $I$-bundle which is \emph{essential with respect to $\mathcal{A}$} if (i) $\rho(S \times 0)$ and $\rho(S \times 1)$ are essential subsurfaces of $\partial M \setminus \mathcal{A}$, (ii) $\mathcal{A} \subset \rho(\partial S \times I)$ and (iii) for any $x \in S$, $\rho(x \times I)$ cannot be homotoped rel boundary into $\overline{\partial M \setminus \mathcal{A}}$.  The following is Lemma 2.1 of \cite{li}: 

\begin{lemma} \label{product} Let $(M, \mathcal{A} )$ be as above.  Then there is a maximal $I$-bundle region $J$ in $M$ such that every product disk with respect to $\mathcal{A} $ can be properly homotoped in $M$  into $J$. 
\end{lemma}

The proof is by construction; successively adding regular neighborhoods of product disks, and then filling in balls with $I$-bundles.  The result is an $I$-bundle which is essential with respect to $\mathcal{A}$.  In particular, fibers of $J$ cannot be homotoped in $M$ rel boundary into $\overline{\partial M \setminus \mathcal{A}} $.  We will call the resulting $I$-bundle a \emph{maximal product disk $I$-bundle}.

Let $X$ be an irreducible manifold with torus boundary components, and let $S \subset X$ be a two-sided, incompressible surface with boundary. Denote the regular neighborhood of $S$ by $N(S)$.  Then cut $X$ along $S$ to obtain a manifold $ X \setminus N(S)$ with boundary and a collection of annuli $\mathcal{A} $.  The union of the collection of annuli is the set of boundary tori in $X$ that meet $S$ minus their intersection with $N(S)$.  Then $X \setminus N(S)$ and $\mathcal{A} $ satisfy (1) and (2) above.  

\begin{definition} \label{essrect} An  \emph{essential rectangle for $(X,S)$} of length $n$  is a map $i: I \times [0,n] \rightarrow X$ such that each $i|_{I \times [k,k+1]}, k \in [0, n-1]$ is a product disk for $X \setminus N(S)$ with respect to $\mathcal{A} $ Furthermore, $i$ is not properly homotopic (fixing $i(I \times 0)$ and $i(I \times n)$) to an essential rectangle of shorter length. 
\end{definition} 

For the purposes of this map, we identify $X$ with the quotient of $X$
 obtained by identifying the $I$-fibers of $N(S)$ to points, so that the map is continuous.  If $S$ is a fiber or a semi-fiber, then there are essential rectangles for $(X,S)$ of arbitrarily high length.   This is not the case for surfaces that are not virtual fibers. 

\begin{lemma}
\label{productdisk}
Let $X$ be an atoroidal, acylindrical, irreducible, compact 3-manifold with torus boundary components and $S$ a two-sided, incompressible surface in $X$.  Suppose that $S$ is not a virtual fiber.  Then there exists a number $P(S) \in \mathbb{N}$ such that the length of any essential rectangle for $(X, S)$ is less than $P(S)$. 
\end{lemma}

The proof follows \cite{li}, although there the statement is only for manifolds with one torus boundary component.  We give the proof in Appendix \ref{onlyone}  for the convenience of the reader.

\section{Normalization} \label{normalization}
In this section we give a procedure to normalize a surface $S$ that satisfies the hypotheses of Theorem \ref{main}.  In Section \ref{finite} below, we show that if this procedure does not normalize $S$ in a finite number of steps, then there exist arbitrarily long essential rectangles for $(M,S)$, which is a contradiction to the assumption that $S$ is not a virtual fiber. 

\begin{definition}Let $N$ be pseudo-manifold. We say that an isotopy of  the image of a surface $\Sigma$ in $N$ \emph{respects the end structure} if the restrictions of the isotopy to $\Sigma \cap N_0$ and $\Sigma \cap \mathcal{E}$ are proper and the isotopy restricted to the 0-skeleton of $N$ is fixed throughout the isotopy. 
\end{definition}

Let $\mathcal{T}: M \rightarrow N$ be a triangulation and $\Sigma$ the image of a surface in $N$ such that  $\phi_{\mathcal{T}}(\Sigma \cap N_0)$ is properly isotopic to a surface $S$ in $M$. Isotopies that respect the end structure will preserve the property that $\phi_{\mathcal{T}}(\Sigma \cap N_0)$ is properly isotopic to $S$. 
\begin{description}

 \item{Step 0:} Initial normalization.

Let $S$ be  a properly embedded incompressible, surface with boundary in a manifold $M$.  Let  $\mathcal{T}: M \rightarrow N$ be an ideal triangulation of $M$ with essential edges. (Recall Definition \ref{ideal} and subsequent discussion.) We assume that $M$, $S$ and $\mathcal{T}$ satisfy the hypotheses of Theorem \ref{main}.  Label the ends of $(M, \mathcal{T})$ by $\mathcal{E}_i$,  call their union $\mathcal{E}$ and denote $\overline{N \setminus \mathcal{E}}$ by $N_0$.  Fix a homeomorphism $\phi_{\mathcal{T}}: N_0 \rightarrow M$.  Let $\Sigma_0 = \phi_{\mathcal{T}}^{-1}(S)$. Thus we have identified $S$ in $M$ with a properly embedded surface $\Sigma_0$ in $N_0$, and $N_0$ is homeomorphic to $M$.   Recall that the universal cover $\widetilde{\mathcal{E}_i}$ of $\mathcal{E}_i$ is $\mathbb{R}^2 \times \mathbb{R}_+$. Let $T_i$ be the boundary component of $N \setminus \mathcal{E}_i$, such that the universal cover $\widetilde T_i$ of each $T_i$ is identified with the plane $z =0$.  The pseudo-manifold $N$ naturally induces a triangulation of each $T_i$.  Since we have no vertices of this triangulation of valence $< 3$, we may assume that this triangulation lifts to a triangulation of $\mathbb{R}^2$ in $\mathbb{R}^2 \times \mathbb{R}_+$ where the edges are linear in the product structure.  In particular, we have no degenerate triangles. A typical line will be normal in this triangulation of $\mathbb{R}^2$. 

Since $\Sigma_0$ is incompressible in $N_0$, the intersection of $\Sigma_0$ with each $T_i$ is either empty or a union of essential curves of slope $p_i/q_i$. We adjust the surface $\Sigma_0$ by an isotopy which respects the end structure so that (1) the surface is in general position with respect to the cell decomposition of $N_0$ induced by $N$, and (2) the boundary curves are normal curves (with respect to the triangulation induced by $N$) of slope $p_i/q_i$ in each $T_i$.  We may assume, after performing another isotopy, that (3) these boundary curves are ``straight". By this we mean that they are covered by lines in the fixed Euclidean product structure on $\widetilde T_i$ which are normal in the lifted triangulation and of slope $p_i/q_i$.   Then by linear changes of coordinates, the boundary curves in each $T_i$ are covered by lines of the form $x=c$ and the groups of deck transformations for the coverings $\widetilde T_i \rightarrow T_i$ are generated by $(x,y) \mapsto (x+1, y)$ and $(x, y) \mapsto (x, y+1)$.

 \item{Step 1:}  Spin $S$. 
 
 Then we extend the lines $x=c$ to slanted half-planes given by $x= z+c $ in $\widetilde{\mathcal{E}_i} =  \mathbb{R}^2 \times \mathbb{R}_+$. Note that we have made a choice of direction to slant the plane.  This choice does not matter, except that it should be consistent for the lifts of the components of $\Sigma \cap \mathcal{E}_i$.  Then the image of each collection of slanted half-planes in $\mathcal{E}_i$ is a union of annuli in the end $\mathcal{E}_i$ that spin once around the cusp as they go out one level in the product structure. Let $\mathtt{T}$ be a tetrahedron in $N$. Then $\Sigma \cap \mathcal{E}_i\cap \mathtt{T}$ is normal in $ \mathcal{E}_i\cap \mathtt{T}$ and either contains an infinite collection of normal triangles or is empty. For each end $\mathcal{E}_i$ which has non-empty intersection with $\Sigma$, we add the vertex in $N$ associated with $\mathcal{E}_i$ to $\Sigma$, and we continue to call the resulting surface $\Sigma$.  Thus we may regard $\Sigma$ as the image of a surface with boundary, where the boundary components map to the vertices of $N$.   Do this for each component of $\Sigma \cap \mathcal{E}$ in  $\mathcal{E}$, and we call the image \emph{spun in $(M, \mathcal{T})$}, or equivalently, \emph{spun in $N$}.  
 
 \end{description}
 Note that now the surface $\Sigma$ is spun, but not spunnormal, and $\phi_{\mathcal{T}}(\Sigma \cap N_0)$ is properly isotopic to $S$.  
 
  \begin{description}
  \item{Step 2:} Remove circle intersections with faces and cancel pairs of intersections with edges.  
  
The goal is to isotope $\Sigma$ respecting the end structure so that $\Sigma$ has no circle intersections with the faces of $N$, and so that each disk of $\Sigma$ intersect a tetrahedron $\mathtt{T}$ has boundary which intersects each edge of $\mathtt{T}$ at most once. Then $\Sigma$ will be normal or spunnormal in $N$.   Here we describe the general procedure.  In Section \ref{finite} below we give more details and prove that this procedure terminates in a finite number of steps. 

It follows from the construction of  $\Sigma \cap \mathcal{E}$  that there are no circle intersections of the faces of $N$ and $\Sigma \cap \mathcal{E}$.   Therefore there are finitely many circle intersections of $\Sigma$ with the faces of $\mathcal{T}$.  Since $S$ is incompressible in $M$,  $\Sigma \setminus N^{(0)}$ is incompressible in $N \setminus N^{(0)}$. We remove these circle intersections   by incompressibility of $\Sigma$,  (and irreducibility of $N \setminus N^{(0)}$) starting with an innermost circle in a face.   Clearly this can be done by an isotopy respecting the end structure.
  
 Let $\mathtt{T}$ be a tetrahedron of $N$. A component of $\Sigma \cap \mathtt{T}$  may have boundary which intersects an edge of $\mathtt{T}$ more than once, say in two points $x$ and $y$.  If this happens, then there is an arc on $\Sigma \cap \mathtt{T}$ that connects $x$ and $y$.  Starting with an innermost such pair along the edge, we can remove the pair of intersections $x$ and $y$, by isotoping across the disk in $\mathtt{T} \setminus \Sigma$ bounded by the arc on $\Sigma \cap \mathtt{T}$ and the arc on the edge.  We can perform such an isotopy without creating any new intersections of $\Sigma$ and the one-skeleton of $N$, but we may have created new circles of intersection with the faces. We then remove these circles of intersection by incompressibility.    In the case that $\Sigma$ is the image of a closed surface, the number of intersections of $\Sigma$ and the edges of $\mathcal{T}$ is finite, and the process described above strictly decreases the number of such intersections.  Therefore the process must end in a finite number of steps. This is Kneser's argument.  Below, we argue that even when the intersection of $\Sigma$ and the the edges of $\mathcal{T}$ is infinite, if $S$ is not a virtual fiber, then the process ends in a finite number of steps.   
   
 \end{description} 
 
 \section{Finiteness} \label{finite} 
 
 \begin{definition} Let $N$ be a pseudo-manifold  and $\Sigma$ a surface in $N$.  A \emph{canceling pair} for $\Sigma$ in $N$ is a pair $(x,y)$ of points on the interior of an edge $e$ of $N$ such that the arc $[x,y]$ on $e$ is isotopic in $N \setminus N^{(0)}$ fixing endpoints to an arc on the surface $\Sigma \setminus N^{(0)}$.  A \emph{canceling disk for the pair $(x,y)$} is an immersion $i: D \rightarrow N$ such that $\partial D$ is the union of two curves $a$ and $b$ with $i(a) = \alpha$, $i(b) = \beta$, where $\alpha$ is  a curve from $x$ to $y$ on $\Sigma$ and $\beta$ is $[x,y] $ on $e$. \end{definition} 
  
The process of isotoping $S$ across a canceling disk for $(x,y)$  to remove the pair of intersections and then removing the (finitely many) circle intersections that may arise is called \emph{canceling the pair $(x,y)$}, as implicit in Step 2 above.  Such an isotopy can be performed that respects the end structure.  Furthermore, we can perform each isotopy so that no new canceling pairs are introduced. Given this, and the fact that the manifold $N_0$ is canonically isomorphic to $M$, Theorem \ref{cancel} below combined with the procedure in Section 3 implies Theorem \ref{main}. 

Generically, a canceling disk is transverse to $\Sigma$ and $\lbrace T_i \rbrace$, and $i^{-1}(\Sigma) \cap i^{-1}(\lbrace T_i \rbrace)$ is a collection of points in $D$. We call the number of such intersections (counted without sign) the \emph{complexity} of $i$ with respect to $\Sigma$ and $\lbrace T_i \rbrace$.  We will show that if there a large number of canceling pairs on a single edge, then any canceling disk for some pair has a high complexity.  This will yield a contradiction to Lemma \ref{productdisk}.

 \begin{theorem} \label{cancel}
 Let $N$ be a pseudo-manifold with essential edges such that $N_0 = \overline{N \setminus \mathcal{E}}$ is an atoroidal, acylindrical, irreducible, compact  3-manifold with torus boundary components. Let $\Sigma$ be surface which is spun in $N$ such that $\Sigma \cap N_0$ is properly embedded, two-sided and incompressible in $N_0$.  If $\Sigma \cap N_0$ is not a virtual fiber of $N_0$, then $\Sigma$ has a finite number of canceling pairs in $N$.\end{theorem} 
 
\begin{proof} We claim that a canceling pair $(x,y)$ in $e \cap \mathcal{E}$ must have $x$ and $y$ in different components of $e \cap \mathcal{E}$.    Suppose that we have a canceling pair in some component of $e \cap \mathcal{E}$, and hence in the same component of $e \cap \mathcal{E}_i$ for some $i$.  Then if the canceling disk $D$ for this pair is contained in $\mathcal{E}_i$, it will lift to the universal cover of $\widetilde{\mathcal{E}_i}$ of $\mathcal{E}_i$  This is a contradiction since a lift of $\Sigma$ is a tilted half-plane and a lift of $e$ is a vertical line, which intersect at most once.  Hence the canceling disk is not contained in $\mathcal{E}_i$ and must intersect $T_i$.   We can remove any circle intersections of $T_i$ and $D$ by incompressibility of $T_i$.  An outermost arc intersection of $T_i$ and $D$ either can be removed or defines a boundary compressing disk of $\Sigma_0 = \Sigma \cap N_0$, which is a contradiction.  This proves the claim. 

As $N_0$ is compact, there are at most finitely many canceling pairs with one or two points in $N_0$, and we cancel these, starting with an innermost pair, by an isotopy which respects the end structure.  This can be done  without introducing any new intersections of $\Sigma$ and $N^{(1)}$.   

Now we have only canceling pairs $(x,y)$ with $x$ in one component of $e \cap \mathcal{E}$ and $y$ in another. Consider some edge $e$ of $N$. We will show that there are less than $P(\Sigma_0)$ canceling pairs on $e$, where $P(\Sigma_0)$ is the maximal length of an essential rectangle as defined in Section 2.   Label the canceling pairs on $e$  $\lbrace (x_j, y_j) \rbrace$, where $(x_0,y_0)$ is an innermost pair and the others are labeled successively.  

Consider the pair $(x_r, y_r)$ and a canceling disk $i: D \rightarrow N$ for the pair $(x_r, y_r)$ with image $D_r$. 
We claim that the minimum complexity of $i$ is $2r+2$.

To prove the claim, we will restrict ourselves to such mappings that have no circles of  $i^{-1}(\Sigma)$ or $i^{-1}(\lbrace T_i \rbrace)$ since removing these  (by incompressibility of $\Sigma \setminus N^{(0)}$ and $\lbrace T_i \rbrace$) will not raise the complexity.  

First note that there is disk with exactly $r +1$ arcs of $i^{-1}(\Sigma)$, one for each pair $\lbrace (x_j, y_j) \rbrace$, which we call $a_j$.  ($a_r = a$.)  Indeed, there are at least $r+1$ arcs of $i^{-1}(\Sigma)$ since there are $r+1$ pairs of intersections on $e$. Any additional arc would give a canceling pair with both intersection points in one component of $e \cap \mathcal{E}$, or a canceling pair with at least one point in $e \cap N_0$. The former doesn't exist by the argument above and the latter we have already removed  by assumption.  

In a minimal complexity disk, there are exactly two arcs of $i^{-1}(\lbrace T_i \rbrace)$.  Indeed, the $T_i$ are normal vertex-linking tori, and hence intersect each edge of an ideal triangulation twice. Since the edges are essential, these intersection points are not connected by an arc on any $T_i$ which bounds a disk with $e \cap N_0$. Therefore there are  at least two arcs of $i^{-1}(\lbrace T_i \rbrace)$ in $D$.  $\Sigma \cap N_0$ and $(\Sigma \setminus N^{(0)}) \cap \mathcal{E}$ are both $\partial$-incompressible with respect to $\lbrace T_i \rbrace$.  Therefore, we can isotope $D$ to remove any additional arcs of  $i^{-1}(\lbrace T_i \rbrace)$ (starting with an outermost such arc) and reduce the number of intersections of $i^{-1}(\Sigma)$ and $i^{-1}(\lbrace T_i \rbrace)$. Thus in a minimal complexity disk there are exactly two arcs of 
$i^{-1}(\lbrace T_i \rbrace)$, which we call $t_1$ and $t_2$. 

Each $t_i$ intersects each $a_j$ at most once in a minimal complexity disk. Indeed, consider the canceling subdisk of $D$, $D_j$, for the canceling pair $(x_j,y_j)$.  If some $t_i$ intersects $a_j$ in more than one point, we can reduce the number of arcs of $i^{-1}(\lbrace T_i \rbrace)$ for $D_j$ by the argument above, reducing the complexity of $D$.  Also, each $t_i$ intersects each $a_j$ at least once since  $a_j$ connects the components of $e \cap \mathcal{E}$.  This proves the claim that the minimal complexity canceling disk for $(x_r,y_r)$ has complexity $2r+2$. 
\vskip .2 in

Call the segment of $a_j$ connecting $t_1$ and $t_2$ $a_j(t_1,t_2)$.   Similarly, we call the segment of $t_i$ connecting $a_k$ and $a_l$  $t_i(a_k, a_l)$.  It intersects exactly $k-l-1$ other arcs $a_j$.  

Let $i: D \rightarrow N$ define a canceling disk of minimal complexity for the pair $(x_r,y_r)$ as above. We claim that the restriction of $i$ to a disk in $D$ bounded by the four arcs $t_1(a_0, a_r)$, $a_{r}(t_1,t_2)$, $t_2(a_r,a_0)$ and $a_0(t_2,t_1)$ is an essential rectangle of length $r$ for $(N_0,\Sigma_0)$ as in Definition \ref{essrect}.  Each product disk is the restriction of $i$ to the disk $D_j \subset D$ bounded by  $t_1(a_j, a_{j+1})$, $a_{j+1}(t_1,t_2)$, $t_2(a_{j+1},a_j)$ and $a_j(t_2,t_1)$.  The $t_i$ are the vertical boundary components and the $a_j$ are the horizontal boundary components.  We must show (1) that the images of $t_1(a_j, a_{j+1})$ and $t_2(a_{j+1}, a_j)$ are essential arcs in $\mathcal{A}$, and (2) that the images of $a_{j+1}(t_1,t_2)$ and $a_j(t_2,t_1)$ cannot be homotoped rel boundary into $\mathcal{A}$. 
 
 For (1), Suppose that the image of $t_1(a_j, a_{j+1})$ can be homotoped into $\Sigma$.  Then, since  $\Sigma_0$ and $\Sigma \cap \mathcal{E}$ are both  $\partial$-incompressible with respect to $\lbrace T_i \rbrace$, we can find a disk  $D$ with fewer intersections of $i^{-1}(\Sigma)$ and $i^{-1}(\lbrace T_i \rbrace)$.  Similarly for (2), if the image of $a_{j+1}(t_1,t_2)$ is homotopic into $\mathcal{A}$, we can find a disk with lower complexity than $D$.  If the restriction of $i$ is properly homotopic (fixing $a_{r}(t_1,t_2)$ and $a_0(t_2,t_1)$) to an essential rectangle of shorter length, then there is a disk with with lower complexity, defined by this essential rectangle. 
 
Thus $i$ restricts to an essential rectangle of length $r$ for $(N_0, \Sigma_0)$.  By Lemma \ref{productdisk}, $r \leq P(\Sigma_0)$.  Therefore, there are less than $P(\Sigma_0)$  canceling pairs on an edge $e$.  Since there are finitely many edges, this proves Theorem \ref{cancel}. 
 \end{proof}
 \noindent \emph{Proof of Theorem \ref{main}}  Let $M$, $\mathcal{T}:M \rightarrow N$, and $S$ be as in the statement of Theorem \ref{main} with a fixed homeomorphism $\phi_{\mathcal{T}}: N_0 \rightarrow M$.  Then $\Sigma_0 = \phi_{\mathcal{T}}^{-1}(S)$ is a properly embedded two-sided incompressible surface in $N_0$ which is not a virtual fiber.  We then spin $\Sigma_0$ and normalize as in Section \ref{normalization}.  By Theorem \ref{cancel}, this process terminates.  The resulting surface $\Sigma$ is spunnormal, and $\phi_{\mathcal{T}}(\Sigma \cap N_0)$ is properly isotopic to $S$. \qed
 
 \section{Further directions}  There are several very natural questions. 
 \begin{description} 
 \item (1) Is there a ``good" form for fibers in cusped hyperbolic manifolds?  For example,  can every fiber in an ideal triangulation with essential edges be isotoped to either be spunnormal or a union of faces in the triangulation?  The later possibility happens in a layered triangulation. This might allow one to efficiently recognize when a triangulated cusped manifold was fibered. See \cite{CooperTillmann} for an algorithm to determine the fibered faces of the Thurston norm ball for a closed 3-manifold using normal surface theory. 
 
 \item (2) Is there a way to recognize when a surface in spunnormal form is an incompressible surface? 
 \item (3) Is there an analogous theory for immersed surfaces?  In particular, given answers to questions (1) and (2), could one test for an immersed incompressible surface? Could one detect when it was a fiber?  This seems very difficult in general.

 \end{description} 
 An algorithm to recognize when an embedded surface is incompressible and some extensions of the theory to immersed surfaces are referenced to in future work in \cite{Rubspun}.
\bibliographystyle{alpha}
\bibliography{rewritten.bib}

\appendix
\section{Proof of Lemma 2.4}
\label{onlyone}

\vskip .2 in 
{\bf Lemma \ref{productdisk}}
{\it Let $X$ be an atoroidal, acylindrical, irreducible, compact  3-manifold with torus boundary components and $S$ a two-sided, incompressible surface in $X$.  Suppose that $S$ is not a virtual fiber.  Then there exists a number $P(S) \in \mathbb{N}$ such that the length of any essential rectangle for $(X, S)$ is less than $P(S)$. }
\vskip .2 in
The proof follows \cite{li}. 
By taking two parallel copies of $S$ if necessary,  $S$ is separating. Without loss of generality,  $X \setminus S$ has exactly two components, $M_1$ and $M_2$.  

Let $\mathcal{A} _1$ and $\mathcal{A} _2$ be the annuli in $M_1$ and $M_2$, respectively, that are the result of cutting $\partial X$ along the curves $\partial S$.  Let $J_1$ be the maximal product disk $I$-bundle for $M_1$ and $\mathcal{A} _1$ and $J_2$ the maximal product disk $I$-bundle for $M_2$ and $\mathcal{A} _2$, as in Lemma \ref{product}.    Let $S_i = \partial M_i \setminus \mathcal{A} _i$ be the copy of $S$ on each component $M_i$, and $C_i = \partial J_i \cap S_i$.  Since $X$ is not a virtual fiber, at least one of the $C_i$, say $C_1$, is not equal to all of $S_1$. There is an involution $\tau_i, i \in \lbrace 1,2 \rbrace$, that takes  $C_i$ to itself, by switching the endpoints of each fiber in the fiber bundles $J_i$.  Denote the gluing map by $\phi: S_1 \rightarrow S_2$.

Let $R_1$ and $R_2$ be subsurfaces of a surface $S$.   We say that $R_1$ and $R_2$ are equivalent ($R_1 \sim R_2$)  if $R_1$ is isotopic to $R_2$ after adding disk components of the complements $S \setminus R_i$ and then removing any disk components from the resulting surfaces.    By \cite{li}[Prop. 2.2],  if $R_1$ and $R_2$ are subsurfaces of $S$ with $\partial S \subset R_1 \cap R_2$, then there are subsurfaces $R'_1 \sim R_1$ and $R'_2 \sim R_2$ such that if a non-trivial curve can be homotoped into $R_1$ and $R_2$, it can be homotoped into $R'_1 \cap R'_2$.    We denote the set of surfaces equivalent to $R$ by $[R]$, and do not distinguish between $[R]$ and a properly chosen element of $[R]$, so that $[R_1] \cap [R_2]$ denotes $[R'_1 \cap R'_2]$. 

If $C_i' \subset [C_i]$, then we can isotope $J_i$ so $J_i \cap S_i =C_i'$, and define $\tau_i$ coherently.  Then  $\tau_1( [\phi^{-1}(\tau_2([\phi(C_1)] \cap [C_2]))] \cap [C_1])$  is the part of $[C_1]$ that contains both horizontal boundaries of  essential rectangles of length 2 that start in $M_2$.
Let $F_1 = C_1$,  $E_k= \tau_1 ( 
[\phi^{-1}(\tau_2([\phi(F_k)]\cap [C_2]))] \cap [C_1])$ and $F_{k+1} = F_k \cap E_k$.   Note that a regular neighborhood of  $\partial \mathcal{A} _1 \subset S_1$ is contained in every $F_k$.  Call this neighborhood $\mathcal{A}_{S_1}$ The main point of the proof is that unless $[F_k] = [\mathcal{A}_{S_1}]$, $[F_{k+1}] \subsetneq [F_k]. $

Suppose $[F_{k+1}] = [F_k].$  Then 
\begin{description} 
\item{(1)}  $[\phi(F_k)] = [\phi(F_k)] \cap [C_2]$,
\item{(2)} $[\phi^{-1}(\tau_2([\phi(F_k)]\cap [C_2]))] = [\phi^{-1}(\tau_2([\phi(F_k)]\cap [C_2]))] \cap [C_1]$, and 
\item{(3)} $[E_k] =[F_k]$. 
\end{description}

Since we are assuming $C_1 \neq S_1$, there is a non-trivial boundary curve of $[F_k]$ unless $[F_k] = [\mathcal{A} _{S_1}]$.  Call it $\gamma_0$.  Then $\phi(\gamma_0)$ is also a boundary component of $[\phi(F_k)] \cap [C_2]$, by (1), and $\tau_1(\phi^{-1}(\tau_2(\phi(\gamma))))$ is a boundary component of $[E_k] = [F_k]$ by (2) and (3).  Call this new curve $\gamma_1$.   Then, by the same reasoning, $\tau_1(\phi^{-1}(\tau_2(\phi(\gamma_1))))$ is also a boundary curve of $F_k$, and we define $\gamma_i$ this way for all $i$.   Suppose $\gamma_i$ is isotopic to $\gamma_{i+1}$.  Then $\tau_2(\phi(\gamma_i)) \cup \phi(\gamma_i)$ and $\phi^{-1}(\tau_2(\phi(\gamma_i))) \cup \gamma_{i+1}$  bound two annuli or M\"obius strips which glue together to form an torus or a Klein bottle. Otherwise, each $\gamma_i$ cobounds an annulus with $\gamma_{i+1}$ or one of $\gamma_i$ or $\gamma_{i+1}$ bounds a M\"obius strip.   Since $[F_k]$ only has finitely many boundary curves, eventually some $\gamma_i$ is the same boundary curve as $\gamma_j$ for some $j  \neq i$.  This will yield two annuli or two M\"obius strips which glue together. Thus there is an immersion of a torus (or Klein bottle) into $M$. We denote both the torus (or Klein bottle) and its image by $T$.

We claim that $T$ is incompressible and not homotopic into the boundary in $X = M_1 \cup M_2$.  Suppose that there is a compressing disk $f:D \rightarrow X$ with $f(\partial D)$ a non-trivial curve on $T$.  Choose $f$ so that $f(D)$ is transverse to $S$ and $f^{-1}(S \cap f(D))$ is a collection of arcs of minimal cardinality.  If $f(\partial D)$ is isotopic to some $\gamma_i$, then this contradicts the incompressibility of $S$, since $\gamma_i$ is a non-trivial curve of $F_k \subset S$.  Therefore, $f(\partial D)$ must intersect each $\gamma_i$ on $T$ and $f(\partial D) \cap M_i$ is a collection of sub-arcs of $f(\partial D)$  that can be homotoped rel boundary to vertical arcs on  $T \cap J_i, i = \lbrace 1,2 \rbrace$. The pullback of $f(D) \cap S$ to $D$ is a collection of arcs in $D$ that connect the pullbacks of these sub-arcs in $\partial D$.   An outermost such arc in $D$ cuts off a subdisk $D'$ of $D$. $f(D')$ is bounded by an arc of $S$ union a vertical  arc of $T \cap J_1$ or $T \cap J_2$.  This is a contradiction, since the vertical arcs cannot be homotoped rel boundary into $S$.  Therefore, $T$ is incompressible in $X$.  

Suppose that $T$ is homotopic into the boundary.  Then $T$ is an immersed torus and there is a map $\phi: T^2 \times I \rightarrow X$ such that $\phi(T^2 \times 0) = T$ and $\phi(T^2 \times 1)$ is a component $Y$ of $\partial X$. If $Y \cap S = \emptyset$ then, since $S$ is separating,  there is a vertical arc of $J_1$ or $J_2$ that is homotopic into $S$ rel boundary, which is a contradiction.  Suppose $Y$ contains some components of $\partial S$.  $T \cap S$ is a union of non-trivial curves on $S$ which are not homotopic to each other on $S$. Then the pullback of $\phi(T^2 \times I) \cap S$ contains annuli $a_i$ with the property that one boundary curve of each  $a_i$ maps to a component of $Y \cap S$ and the other boundary curve maps to a component of $T \cap S$, for all components of $T \cap S$.  Consider some such annulus $a$ where one boundary curve maps to  $\gamma$.  Then $\phi \vert_a$ defines a homotopy of $\gamma$ into $\partial S$ in $S$, and hence $\gamma$ is isotopic into $\partial S$ by \cite{Ep}.  This is a contradiction since $\gamma$ was not boundary parallel.   Thus $T$ is an incompressible torus (or Klein bottle) that is not homotopic into the boundary, contradicting the fact that $X$ is atoroidal.

\vskip .5 in 
\noindent
Genevieve Walsh\\  Mathematics  \\ Tufts University \\ Medford, MA 02155\\ genevieve.walsh@tufts.edu

\end{document}